\documentclass[11pt]{article}

\usepackage{amsfonts}
\usepackage{amssymb,amsmath}
\usepackage{latexsym}
\usepackage{epsfig}

\newcommand     {\RR}{\mathbb{R}}
\newcommand     {\PP}{\mathbb{P}}
\newcommand     {\EE}{\mathbb{E}}
\newcommand     {\Un} {{\bf{1}}}
\newcommand {\lbd}      {\lambda}

\newcommand     {\vareps}       {\varepsilon}

\newtheorem     {thm}           {Theorem}[section]

\begin{document}


\title{Adaptive dynamics in logistic branching populations}


\author{Nicolas Champagnat\footnote{Project-team TOSCA, INRIA Sophia
    Antipolis\newline e-mail: Nicolas.Champagnat@sophia.inria.fr}, Amaury
  Lambert\footnote{Unit of Mathematical Evolutionary Biology,
    Laboratoire d'{\'E}cologie et {\'E}volution, Universit{\'e} Pierre
    et Marie Curie-Paris 6 and \'{E}cole Normale Sup{\'e}rieure,
    Paris\newline e-mail: amaury.lambert@ens.fr}}

\maketitle

\textit{MSC 2000 subject.} Primary 60J80; secondary 60J25,
  60J70, 60J75, 60J85, 60K35, 92D10, 92D15, 92D25, 92D40.\\


\section{Introduction}
The recent biological theory of adaptive dynamics~[6,9] proposes a
description of the long term evolution of an asexual population by
putting emphasis on the ecological interactions between individuals,
in contrast with classical population genetics models which focus on
the genetic structure of the population. The basic models are
individual-based models in which the population dynamics is precisely
described and includes birth, death, competition and mutation. The
basic idea of the theory of adaptive dynamics is to try to get
insights about the interplay between ecology and evolution by studying
the invasion of a single mutant type appearing in a simplified
(monotype stable) resident population.  The evolution of the
population can then be described as a sequence of mutant invasions in
the population. If the resident type goes extinct when the mutant type
invades (we say that the mutant type \emph{fixates}), the evolution is
described by the so-called `trait substitution sequence' (TSS)~[10].
This appraoch has revealed powerful to predict the qualitative
behaviour of complicated evolutionary dynamics. In particular, it
allows to determine the (local) direction of evolution in the space of
phenotypic traits (or simply \emph{traits}) from the individual
ecological parameters, and to predict and explain the phenomenon of
evolutionary branching~[10], where a population, initially composed of
individuals with traits concentrated around a single trait value,
divides into two (or more) subpopulations concentrated around distinct
trait values that stably coexist because of their ecological
interactions. The description of this phenomenon is an important
achievement of this theory, as well as the `canonical equation of
adaptive dynamics'~[3], which describes the evolution of the dominant
trait of the population as a deterministic ``hill-climbing'' process
on a fitness landscape which depends on the current state of the
population (see~(\ref{eqn:CEAD}) below).

More formally, as soon as eternal coexistence of two or more types is
not permitted by the model, the evolution can be described by the
sequences $(T_n)_n$ and $(V_n)_n$, where $T_n$ is the $n$-th time
where the population becomes monomorphic (i.e.\ composed of only one
type) and $V_n$ is the surviving type at time $T_n$. The sequence
$(V_n)$ is the above-mentionned TSS. It is possible to prove the
convergence of an individual-based model to the TSS under two
biologically motivated assumptions~[10,1]. First, the \emph{assumption
  of rare mutations} guarantees that, in the timescale of mutations,
the widths of time intervals during which the population is
polymorphic vanish, so that there is only one type surviving at any
time $t$. To prevent the population from rapidly becoming extinct in
the new timescale, one also has to rescale population sizes, thereby
making the \textit{assumption of large populations}.

Subsequently, the TSS is a Markov jump process on the trait space
whose semigroup is shown~[1] to depend on the \textit{invasion
  fitnesses} (as defined in~[9]) $f(x,y)$, $x,y \in{\cal X}$, where
$f(x,y)$ is defined as the \textit{expected growth rate} of a single
individual of type $y$ --- the mutant --- entering a monomorphic
population of type $x$ `at equilibrium' --- the residents.  Note that
this fitness is not given \emph{a priori}, but derived from the
microscopic model of individual interactions. Because of the
assumption of large population, the sign of this fitness determines
the possibility of invasion of a mutant type: if $f(x,y)<0$, the
mutant type $y$ cannot invade a resident population of type $x$. Thus,
evolution proceeds by \textit{successive invasions of (only)
  advantageous mutant types} replacing the resident one.

The `canonical equation of adaptive dynamics'~[3], which describes the
evolution of a one-dimensional trait $x$ as the solution of the
following ODE, is obtained from the TSS in the limit of \emph{small
  mutations}:
\begin{equation}
\label{eqn:CEAD}
\frac{dx}{dt}=\frac{1}{2}\sigma(x)^2{\mu(x)}
   {\bar{n}(x)}\frac{\partial}{\partial y}{f(x,x)},
\end{equation}
where $\sigma(x)^2$ stands for the (rescaled) variance of the mutation
step law, $\bar{n}(x)$ for the equilibrium size of a pure $x$-type
population, and $f(x,y)$ for the invasion fitness mentioned above.
Note how only advantageous types get fixed (the trait follows the
fitness gradient) and how the fitness landscape $y\mapsto f(x,y)$
depends on the current state $x$ of the population.

However, it is well-known that slightly deleterious types can be fixed
by chance in finite populations. This phenomenon is known under the
name of \textit{genetic drift}. Depending on the strength of genetic
drift, selection is said to be \emph{strong} (genetic drift has
negligible effects) or \emph{weak}. In the large population asymptotic
from which the TSS of adaptive dynamics is derived, genetic drift has
negligible impact compared to the action of selection.  Therefore, the
fixation of slightly deleterious types cannot be observed. Our goal
here is to include genetic drift in the adaptive dynamics models by
considering \emph{finite populations} under \emph{weak selection}. We
continue using the bottom-up approach of adaptive dynamics; that is,
model (macroscopic) evolution from (microscopic) populations. In
particular, we allow the population sizes to fluctuate randomly
through time and we aim to reconstruct a fitness function from the
microscopic parameters.

After the description of the model (Section~2), we derive a new TSS in
the limit of rare mutations (Section~3), from which a limit of small
mutations gives what we call the `canonical diffusion of adaptive
dynamics' (Section~4). The coefficients of this diffusion involve the
first-order derivatives of the fixation probabilities, which are
computed in Section~5 as a linear combination of four fundamental
components associated to \emph{fertility}, \emph{defence},
\emph{aggressiveness} and \emph{isolation}. New numerical results on
the robustness of the population with respect to these fundamental
components are also given, as well as some consequences on the
canonical diffusion of adaptive dynamics in large populations.

\section{The microscopic model}
We will restrict here to logistic interaction. More general models are
considered in~[2]. 

A monotype (binary) logistic branching process (LBP, see~[7]) with
dynamical parameters $(b,c,d)$ is a Markov chain in continuous time
$(X_t;t\ge 0)$ with nonnegative integer values and transition rates
$$
q_{ij}=\;
\left\{
                \begin{array}{cll}
bi &\mbox{if }&j=i+1\\
ci(i-1) &\mbox{if }& j=i-1\\
-i(b+c(i-1)) &\mbox{if }&j=i\\
0 & \mbox{otherwise.} &
                \end{array}
        \right.
$$
The non-linear term $ci(i-1)$ describes competition mortality due to
random encounters between individuals. Other terms correspond to
independent birth events with constant individual rates. This Markov
chain is positive-recurrent and converges in distribution to a r.v.\
$\xi$, where $\xi$ is a Poisson variable of parameter $\theta:=b/c$
conditioned on being nonzero
\begin{equation}
  \label{eq:def-xi}
  \PP(\xi=i)=\frac{e^{-\theta}}{1-e^{-\theta}}
  \:\frac{\theta^i}{i!}
  \qquad i\geq 1.
\end{equation}

We consider in the sequel, a multitype asexual birth and death process
with mutation, generalizing this LBP. At any time $t$, the population
is composed of a finite number $N(t)$ of individuals characterized by
their phenotypic traits $x_1(t),\ldots,x_{N(t)}(t)$ belonging to a
given trait space ${\cal X}$, assumed to be a closed subset of
$\mathbb{R}^k$. The population state at time $t$ is represented by the
counting measure on ${\cal X}$
\begin{equation*}
  \nu_t=\sum_{i=1}^{N(t)}\delta_{x_i(t)}.
\end{equation*}
The population dynamics is governed by the following parameters.
\begin{itemize}
\item $b(x)$ is the rate of birth from an individual of type $x$. The
  function $b$ is assumed to be ${\cal C}^2_b$.
\item $c(x,y)$ is the rate of death of an individual of type $x$ due
  to the competition with another individual of type $y$.  Therefore,
  the total death rate of an individual of type $x$ in a population
  $\nu$ may be written as $\int c(x,y)(\nu(dy)-\delta_x(dy))$. In this
  expression, the Dirac mass at $x$ substracted to the measure $\nu$
  means that the individual does not compete with himself. The
  function $c$ is assumed to be ${\cal C}^2_b$ and bounded away from 0
  on ${\cal X}^2$.
\item $\gamma\mu(x)$ is the probability that a birth from an
  individual with trait $x$ produces a mutant individual, where
  $\mu(x)\in[0,1]$ and where $\gamma\in(0,1]$ is a parameter scaling
  the frequence of mutations. When there is no mutation, the new
  individual inherits the trait of its progenitor. In Section~3, we
  will be interested in the limit of rare mutations
  ($\gamma\rightarrow 0$).
\item $M(x,dh)$ is the law of the trait difference $h=y-x$ between a
  mutant individual with trait $y$ born from an individual with trait
  $x$. We assume that $M(x,dh)$ has 0 expectation (no mutation bias),
  i.e.\ $\int hM(x,dh)=0$, and has a density on $\mathbb{R}^k$ which
  is uniformly bounded in $x\in{\cal X}$ by some function $\bar{M}(h)$
  with finite third-order moment.
\end{itemize}
We will denote the dependence of $\nu_t$ on the parameter $\gamma$
with the notation $\nu^{\gamma}_t$. Observe that such a population
cannot go extinct because the death rate is 0 when there is only one
individual in the population. Since we want to apply a limit of rare
mutations while keeping the population size finite, this is necessary
to prevent the population to become extinct before any mutation occur.

Let $\xi(x)$ be a random variable whose law is the stationary
distribution of a pure $x$-type population with \textbf{no} mutation
($\mu\equiv 0$). This law is given by~(\ref{eq:def-xi}) where $\theta$
is replaced by $\theta(x):=b(x)/c(x,x)$.

The last notation needed concerns a population with initially only two
types $x$ and $y$ and with \textbf{no} mutation.  Then $\nu_t=X_t
\delta_x + Y_t \delta_y$, where $(X_t,Y_t:t\ge 0)$ is a bivariate
Markov chain. For this Markov chain, $\PP ( T<\infty) =1$, where $T$
is the first time where either $X_t$ or $Y_t$ reach 0. We call
\textit{fixation} (of the mutant $y$) the event $\{X_T=0\}$.  The
probability of fixation will be denoted by $u_{n,m}(x,y)$
$$
u_{n,m}(x,y):=\PP(X_T=0 \:\vert\: X_0=n, Y_0=m).
$$

\section{The trait substitution sequence in finite populations}
In this section, we apply the limit of rare mutations
($\gamma\rightarrow 0$) to the process $\nu^\gamma$, in order to
describe the evolution of the population as a TSS in finite
population. This limit requires to rescale time properly, as
$t/\gamma$, to describe the evolution on the mutation timescale.

\begin{thm}
  Fix $x\in{\cal X}$. Assume that
  $\nu^\gamma_0=N_0^\gamma\delta_x$ where
  $\sup_{\gamma\in(0,1)}\EE((N_0^\gamma)^p)<\infty$ for some $p>1$.
  Then, for any $0<t_1<\ldots<t_n$, the $n$-tuple
  $(\nu^\gamma_{t_1/\gamma},\ldots,\nu^\gamma_{t_n/\gamma})$ converges
  in law for the weak topology to
  $(N_{t_1}\delta_{S_{t_1}},\ldots,N_{t_n}\delta_{S_{t_n}})$ where
  \begin{description}
  \item[(1)] $(S_t;t\geq 0)$ is a Markov jump process on ${\cal X}$
    with initial value $S_0=x$ and whose jumping rates $q(x, dh)$ from
    $x$ to $x+h$ are given by
    \begin{equation*}
      q(x,dh) = \beta(x) \chi(x, x+h) M(x,dh),
    \end{equation*}
    where
    $\beta(x)=\mu(x)b(x)\EE(\xi(x))=\mu(x)b(x)\theta(x)/(1-e^{-\theta(x)})$
    and
    \begin{equation}
      \label{eqn:chi}
      \chi(x,y) = \sum_{n\ge 1}
      \frac{n\PP(\xi(x)=n)}{\EE(\xi(x))}u_{n,1}(x,y) 
      = \sum_{n\ge 1} e^{-\theta(x)}\frac{\theta(x)^{n-1}}{(n-1)!}u_{n,1}(x,y) .
    \end{equation}
  \item[(2)] Conditional on
    $(S_{t_1},\ldots,S_{t_n})=(x_1,\ldots,x_n)$, the $N_{t_i}$ are
    independent and respectively distributed as $\xi(x_i)$.
  \end{description}
\end{thm}

Therefore, in the limit of rare mutations, on the mutation timescale,
the population is always monomorphic and the dominant trait of the
population evolves as a jump process over the trait space, where a
jump corresponds to the appearance and fixation of a mutant type.
Moreover, at any time, the population size stationary (i.e.\ has the
stationary distribution corresponding to the dominant trait of the
population). The fixation rate of a mutant is governed by the function
$\chi(x,y)$, which is therefore the random analogue of the traditional
invasion fitness~[9], defined as the probability of invasion of a
mutant type $y$ is a resident population of type $x$ at equilibrium.
Observe that, as usual in adaptive dynamics, the fitness landscape
depends on the current state of the population. Moreover, in contrast
with the classical TSS~[10], from a given monomorphic resident
population, any mutant trait has a positive probability to invade (by
genetic drift).  Therefore, evolution is possible in any direction of
the trait space.  However, a directional selection still exists, as
will appear in the next section.

We refer to~[2] for the proof of this result. However, this
convergence is natural in view of the following interpretation of each
parameter. $\beta(x)$ can be seen as the \emph{mean mutant production
  rate} of a stationary $x$-type population (i.e.\ with size
$\xi(x)$), and $\chi(x,y)$ is the \emph{probability of fixation} of a
single $y$-type mutant entering a pure $x$-type population with
\emph{size-biased stationary size}. The size bias comes from the fact
that the mutant appears \emph{at a birth time} in the stationary
population (since the birth rate is proportional to the population
size, the population size after a birth event in the stationary
population is given by the size-biaised stationary population size).

\section{The canonical diffusion of adaptive dynamics}
Let us assume for simplicity that ${\cal X}=\mathbb{R}^k$. Let
$\sigma(x)$ be the square root matrix of the covariance matrix of
$M(x,\cdot)$. We also need to assume that the matrix $\sigma(x)$ is a
Lipschitz function of $x$.

In order to obtain the equivalent of the canonical equation of
adaptive dynamics in a finite population, we want to apply a limit of
small mutation steps (weak selection) to the TSS $S$. To this aim, we
introduce a parameter $\epsilon>0$ and replace the mutation kernels
$M(x,\cdot)$ by their image by the application $h\mapsto\epsilon h$.
Time also has to be rescaled in order to obtain a non-degenerate
limit. The correct time scaling is $1/\epsilon^2$, which leads to the
following generator for the rescaled TSS $(S^\epsilon_t;t\geq 0)$
\begin{equation*}
  A_\epsilon\varphi(x)=\frac{1}{\epsilon^2}
  \int_{\RR^k}(\varphi(x+\epsilon
  h)-\varphi(x))\beta(x)\chi(x,x+\epsilon h)M(x,dh).
\end{equation*}
Using the assumption that the mutation kernels $M(x,\cdot)$ have 0
expectation, it is elementary to compute the limit of this expression
as $\epsilon\rightarrow 0$ (for sufficiently regular $\varphi$). This
limit, which takes the form of a diffusion generator, explains the
following result (its full proof can be found in~[2]).

\begin{thm}
  If the family $(S^\epsilon_0)_{\epsilon>0}$ has
  bounded first-order moments and converges in law as
  $\epsilon\rightarrow 0$ to a random variable $Z_0$, then the process
  $S^\epsilon$ with initial state
  $S^\epsilon_0$ converges in law for the Skorohod topology on
  $\mathbb{D}(\RR_+,\RR^k)$ to the diffusion process $(Z_t;t\geq 0)$
  with initial state $Z_0$ unique solution to the stochastic
  differential equation
  \begin{equation}
    \label{eq:can-diff}
    dZ_t=\beta(Z_t)\sigma^2(Z_t)\cdot \nabla_2\chi(Z_t,Z_t) dt +\sqrt{\beta(Z_t)\chi(Z_t,Z_t)}\sigma(Z_t)\cdot dB_t
  \end{equation}
  where $\nabla_2\chi$ denotes the gradient w.r.t.\ the second variable
  $y$ of $\chi(x,y)$ and $B$ is a standard $k$-dimensional Brownian motion.
\end{thm}


This result gives the equivalent of the canonical equation of adaptive
dynamics~(\ref{eqn:CEAD}) when the population is finite. It is no
longer a deterministic ODE, but a diffusion process, in which the
genetic drift remains present (in the form of a stochastic diffusion
term), as a consequence of the population finiteness and of the limit
of weak selection. The deterministic drift part of~(\ref{eq:can-diff})
is very similar to the standard canonical equation of adaptive
dynamics (\ref{eqn:CEAD}), and involves in particular the gradient of
the fitness function $\chi$. The process~(\ref{eq:can-diff}) provides
a diffusion model describing the evolution of the dominant trait value
in a population~[8,5], grounded on a precise microscopic
density-dependent modelling of the population dynamics. It also gives
the precise balance between directional selection and genetic drift as
a function of the individual's dynamical parameters.

\section{Fixation probability near neutrality}
The SDE~(\ref{eq:can-diff}) involves the fixation probability
$\chi(x,x)$ and the fitness gradient with respect to the second
variable $\nabla_2\chi(x,x)$. In this section, we explain how these
quantities can be explicitly computed.

We need to compute the derivatives of the fixation probabilities
$u_{n,m}(x,y)$ when $y$ is close to $x$. Recall that the law of the
two-types LBP without mutation $(X,Y)$ used to define $u_{n,m}$ in the
end of Section~2 is characterized by the birth vector $B$ and the
competition matrix $C$
$$
B=\left( 
\begin{array}{c}
b(x)\\ 
b(y)\\ 
\end{array}
\right),\quad
C=\left( 
\begin{array}{cc}
c(x,x) & c(x,y)\\ 
c(y,x) & c(y,y)\\ 
\end{array}
\right).
$$
We will say that the mutant is \emph{neutral} if all individuals are
exchangeable, i.e.\ when $b(y)=b(x)$ and $c(x,y)=c(y,x)=c(y,y)=c(x,x)$
(this holds in particular when $y=x$). As will appear below, using the
notation $(b,c):=(b(x),c(x,x))$, it is natural to focus on deviations
from the neutral case expressed as
$$
B=b\Un+\left( 
\begin{array}{c}
0\\ 
\lbd\\ 
\end{array}
\right),
\quad
C=c\Un-\left( 
\begin{array}{cc}
0& 0\\ 
\delta & \delta\\ 
\end{array}
\right)+\left( 
\begin{array}{cc}
0& \alpha\\ 
0&\alpha  \\ 
\end{array}
\right)-\left( 
\begin{array}{cc}
0& \vareps\\ 
\vareps& 0\\ 
\end{array}
\right)
$$
In words, deviations from the neutral case are a linear combination of
four \textit{fundamental selection coefficients}
$\lambda$, $\delta$, $\alpha$, $\vareps$, that are chosen to
be positive when they confer an advantage to the mutant. It is
convenient to assess these deviations to the neutral case in terms of
\begin{enumerate}
\item \textbf{fertility} ($\lambda$): positive $\lambda$
  means increased mutant birth rate
\item \textbf{defence} capacity ($\delta$): positive $\delta$ means
  reduced competition sensitivity of mutant individuals w.r.t.\ the
  total population size
\item \textbf{aggressiveness} ($\alpha$): positive $\alpha$ means
  raised competition pressure exerted from any mutant individual onto
  the rest of the population
\item \textbf{isolation} ($\varepsilon$): positive $\varepsilon$ means
  lighter cross-competition between the two different types, that
  would lead, if harsher, to a greater probability of exclusion of the
  less abundant one
\end{enumerate}
Under neutrality, an elementary martingale argument shows that the
fixation probability equals the initial mutant frequency $p:=m/(m+n)$.
This implies in particular that
\begin{equation}
  \label{eqn:chi(x,x)}
  \chi(x,x)=\frac{e^{-\theta(x)}-1+\theta(x)}{\theta(x)^2}.  
\end{equation}
The following theorem unveils the dependence of $u$ upon $\lambda$,
$\delta$, $\alpha$, $\vareps$, when they slightly deviate from 0, and
explains why these four selection coefficients provide a natural basis
to decompose the gradient of the fixation probability.

\begin{thm}
  As a function of the multidimensional selection
  coefficient $s = (\lambda, \delta, \alpha, \varepsilon)$, the
  probability $u$ is differentiable, and in a neighbourhood of $s=0$
  (selective neutrality),
  \begin{equation}
    \label{eqn:gradient u}
    u = p + v'.s+o(s),
  \end{equation}
  where the selection gradient $v=(v^\lbd,
  v^\delta,v^\alpha,v^\vareps)$ can be expressed as
  $$
  \begin{array}{lclr}
    v^\iota_{n,m} &=& p \,(1-p)\,\, g^\iota_{n+m}&\iota\not=\varepsilon,\\
    v^\varepsilon_{n,m}& =& p \,(1-p)\,(1-2p)\,\,g^\varepsilon_{n+m}&
  \end{array}
  $$
  where the $g$'s depend solely on the resident's characteristics
  $b,c$, and on the total initial population size $n+m$. They are
  called the invasibility coefficients.
\end{thm}

The invasibility coefficients of a pure resident population are
interesting to study, as they provide insights about how the fixation
probability deviates from $p$ as the selection coefficients of the
mutant deviate from 0. They also provide information about the
\emph{robustness} of the resident population, i.e.\ its resistance to
mutant invasions. In particular, this allows to compare the
sensitivity of the invasion probability in a given monomorphic
resident population with respect to the four fundamental selection
coefficients. In the simplest case where mutations in the parameter
space are isotropic, the biggest invasibility coefficent gives the
direction of the parameter space where a mutant is more likely to
invade. More generally, when there are correlations between mutations
in the parameter space (either because of the phenotypic structure in
the functions $b(\cdot)$ and $c(\cdot,\cdot)$, or because the
covariance matrix of the mutation steps is non-diagonal), the
likeliest direction of evolution in the trait space is given by the
deterministic coefficient of the canonical
diffusion~(\ref{eq:can-diff}), in which the fitness gradient is given
by
\begin{equation}
  \label{eq:grad-chi}
  \nabla_2\chi(x,x)=a_\lambda(x)\nabla
  b(x)-a_\delta(x)\nabla_1c(x,x)+a_\alpha(x)\nabla_2c(x,x),
\end{equation}
where, for $\iota=\lambda,\delta,\alpha$,
\begin{equation*}
  a_\iota(x)=e^{-\theta(x)}\sum_{n=1}^\infty\frac{ng^\iota_{n+1}(x)\theta(x)^{n-1}}
  {(n+1)^2(n-1)!}.
\end{equation*}

\begin{figure}
  \centering
  \input{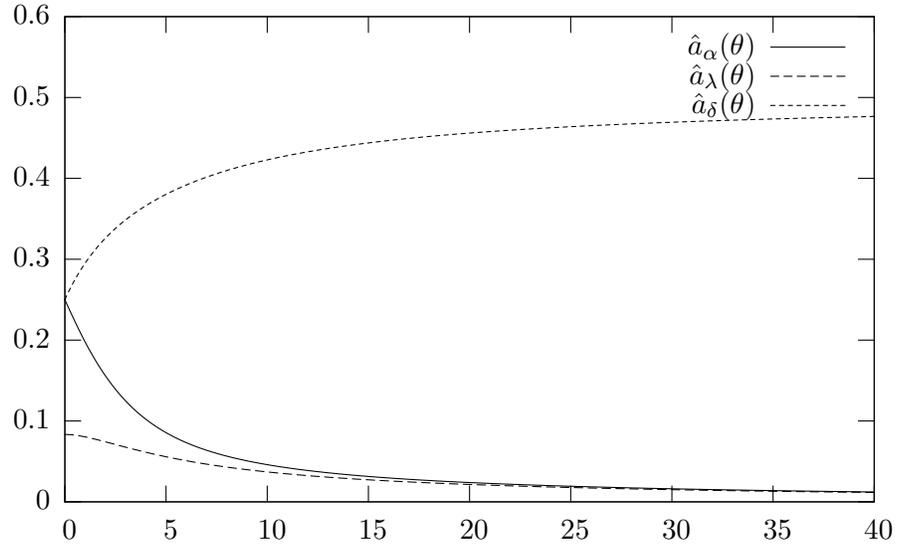}
  \caption{The functions $\hat{a}_\lambda$, $\hat{a}_\delta$ and
    $\hat{a}_\alpha$ as functions of $\theta$.}
  \label{fig:a}
\end{figure}

\begin{figure}
  \centering
  \input{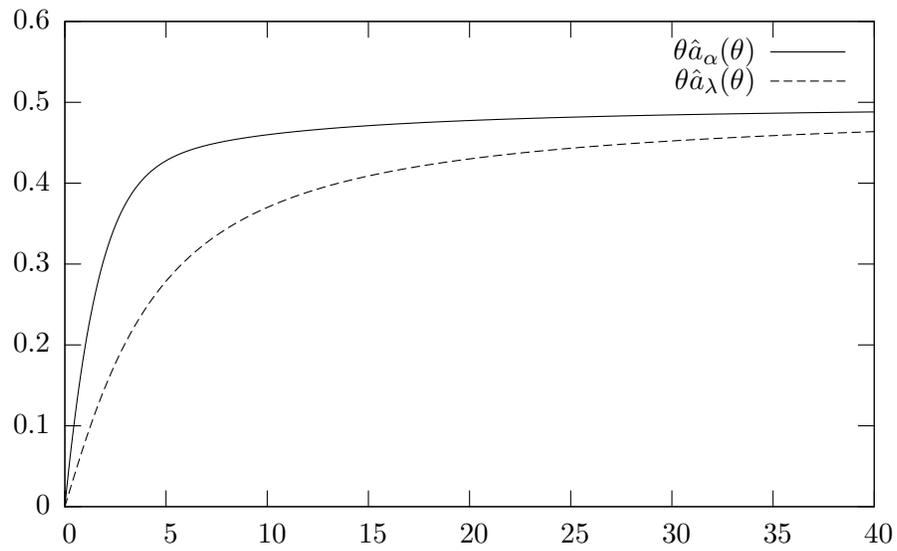}
  \caption{The functions $\theta\mapsto\theta\hat{a}_\lambda(\theta)$
    and $\theta\mapsto\theta\hat{a}_\alpha(\theta)$.}
  \label{fig:a_infty}
\end{figure}

It is possible to obtain explicit expressions for the invasibility
coefficients $g^\iota$ as series. We refer to~[2] for the exact
expressions. In particular, these expressions yield that
$a_\iota(x)=\hat{a}_\iota(\theta(x))/c(x,x)$ for some function
$\hat{a}_\iota$. Moreover, they allow one to compute numerically the
invasibility coefficients, and therefore the quantities
$\hat{a}_\iota$ for $\iota=\lambda,\delta,\alpha$ as functions of the
parameter $\theta(x)$. These numerical results can be used to make
simulations of the canonical diffusion of adaptive dynamics in various
ecological examples. In particular, in contrast with the classical
canonical equation of adaptive dynamics, the presence of a genetic
drift can induce the evolutionary dynamics to drift away from
evolutionary stable strategies, where the fitness gradient is zero.
When the fitness gradient admits several zeros, it can visit several
basins of attraction on various timescales.

This numerical study is a work in progress, that is quite delicate
because the series involved in the computation of $g^\iota$ are slowly
converging, with first terms that grow exponentially fast with
$\theta$. We shall give here our first results. Fig.~\ref{fig:a} shows
the functions $\hat{a}_\iota$ for $\iota=\lambda,\delta,\alpha$.
Several comments can be made from this figure. First, for any
$\theta>0$,
$\hat{a}_\delta(\theta)>\hat{a}_\alpha(\theta)>\hat{a}_\lambda(\theta)$.
This means that, for equal mutation steps in the parameter space, a
mutation is always more advantageous in the direction $\delta$ than in
the direction $\alpha$, which is itself more advantageous than in the
direction $\lambda$. In other words, in a given population, a better
defence capacity is more beneficial than a better aggressiveness,
which is more beneficial than a better fertility.

Moreover, as $\theta$ goes to infinity, these functions have different
asymptotic behaviors. $\hat{a}_\delta$ seems to converges to $1/2$,
whereas $\hat{a}_\lambda(\theta)$ and $\hat{a}_\alpha(\theta)$ are
both equivalent to $1/2\theta$ (see Fig.~\ref{fig:a_infty}). However,
this does not mean that mutations are much more likely to fixate in
the $\delta$ direction, because, when $\theta$ is large, $b=\theta c$
is larger than $c$, and, in~(\ref{eq:grad-chi}), $a_\delta$ and
$a_\alpha$ are multiplied by $\nabla_1 c$ and $\nabla_2 c$
respectively, whereas $a_\lambda$ is multiplied by $\nabla b$. More
formally, to compute the limit of the canonical diffusion when
$\theta$ goes to infinity, one can divide the competition kernel
$c(\cdot,\cdot)$ by a constant $K$ in the microscopic model and then
let $K$ go to infinity in the canonical diffusion. Denoting by
$\chi_K(\cdot,\cdot)$ the fitness function obtained this way and using
the asymptotic behaviors given above, one gets that
\begin{equation*}
  \lim_{K\rightarrow+\infty}\nabla_2\chi_K(x,x)
  =\frac{1}{2b(x)}(\nabla b(x)-\theta(x)\nabla_1 c(x,x)).
\end{equation*}
Moreover, by~(\ref{eqn:chi(x,x)}), $\chi_K(x,x)$ converges to 0 when
$K\rightarrow\infty$. Now, as proved in~[1], with our notation, the
fitness function of the canonical equation of adaptive
dynamics~(\ref{eqn:CEAD}) is given by $f(x,y)=b(y)-c(y,x)\theta(x)$.
Therefore, as $K\rightarrow+\infty$, the canonical diffusion converges
to a deterministic ODE which is precisely the canonical equation of
adaptive dynamics. This gives a new justification of this equation,
and this also allows one to study the fluctuations around the
canonical diffusion when $K$ is large. In particular, this diffusion
equation with small diffusion term enters the framework of
Freidlin-Wentzell's theory~[4], which can be applied to predict the
long time behaviour of the diffusion, and its chain of visit of basins
of attractions when $K$ is large. This kind of information is
biologically very relevant, since it allows one to predict in which
order all the evolutionary stable strategies will be visited by the
population and on which timescale.



\begin{thebibliography}{99}

\bibitem{[1]} N. Champagnat, {\it A microscopic
    interpretation for adaptive dynamics trait substitution sequence
    models}\/, Stoch.\ Proc.\ Appl.\ 116 (2006), 1127--1160.

\bibitem{[2]} N. Champagnat and A. Lambert, {\it Evolution of discrete
    populations and the canonical diffusion of adaptive dynamics}\/,
  to appear in Ann.\ Appl.\ Prob.

\bibitem{[3]} U. Dieckmann and R. Law, {\it The dynamical theory of
    coevolution: A derivation from stochastic ecological processes}\/,
  J.\ Math.\ Biol.\ 34 (1996), 579--612.

\bibitem{[4]} M.I. Freidlin and A.D. Wentzell, Random Perturbations
    of Dynamical Systems. Springer-Verlag, Berlin, 1984.

\bibitem{[5]} T.F. Hansen, {\it Stabilizing selection and the comparative
    analysis of adaptation}\/, Evolution 51 (1997), 1341--1351.

\bibitem{[6]} J. Hofbauer and R. Sigmund, {\it Adaptive dynamics and
    evolutionary stability}\/, Appl.\ Math.\ Letters 3 (1990), 75--79.

\bibitem{[7]} A. Lambert, {\it The branching process with logistic
    growth}\/, Ann.\ Appl.\ Prob.\ 15 (2005), 1506--1535.

\bibitem{[8]} R. Lande, {\it Natural selection and random genetic drift
    in phenotypic evolution}\/, Evolution 30 (1976), 314--334.

\bibitem{[9]} J.A.J. Metz, R.M. Nisbet and S.A.H. Geritz, {\it How should
    we define `fitness' for general ecological scenarios ?}\/, Trends
  in Ecology and Evolution 7 (1992), 198--202.

\bibitem{[10]} J.A.J. Metz, S.A.H. Geritz, G. Mesz{\'e}na, F.A.J. Jacobs
  and J.S. van Heerwaarden, {\it Adaptive Dynamics, a geometrical
    study of the consequences of nearly faithful reproduction}, in:
  Stochastic and Spatial Structures of Dynamical Systems, S.J. van
  Strien and S.M. Verduyn Lunel (eds.), North Holland, Amsterdam,
  1996, 183--231.

\end{thebibliography}
\end{document}